\documentclass{amsart}
\usepackage{graphicx,amsmath,amssymb,amsthm,harvard}
\DeclareGraphicsExtensions{.jpg,.eps,.eps.bb}
\newcommand{\cprime}{$'$}
\newcommand{\new}[1]{{\em #1}}
\newcommand{\lser}{[ [}
\newcommand{\rser}{] ]}
\newcommand{\B}{\mathbb{B}}
\newcommand{\N}{\mathbb{N}}
\newcommand{\R}{\mathbb{R}}
\newcommand{\Z}{\mathbb{Z}}
\newcommand{\sK}{\mathcal{K}}
\newcommand{\sU}{\mathcal{U}}
\newcommand{\sF}{\mathcal{F}}
\newcommand{\sL}{\mathsf{L}}
\newcommand{\sR}{\mathsf{R}}
\newcommand{\rmax}{\mathbb{R}_{\max}}
\newcommand{\rmaxb}{\overline{\R}_{\max}}
\newcommand{\nmaxb}{\overline{\N}_{\max}}
\newcommand{\zero}{\varepsilon}
\newcommand{\unit}{e}
\newcommand{\supp}{\mathop{\text{\Large$\vee$}}}
\newcommand{\inff}{\mathop{\text{\Large$\wedge$}}}
\newcommand{\maxx}{\top}
\newcommand{\minn}{\bot}
\newcommand{\maxxx}{{\top \sK}}
\newcommand{\transpose}{{\text{\tiny$\top$}}}
\newcommand{\leqop}{\mathop{\scriptstyle\stackrel{\mathrm{op}}{\leq}}}
\newcommand{\suppop}{\supp\!\null\op}
\newcommand{\fench}{^{*}}
\newcommand{\op}{^{\mrm{op}}}
\newcommand{\dotop}{\mathop{\scriptstyle\stackrel{\mathrm{op}}{\cdot}}}
\newcommand{\ldual}[1]{{}^-#1}
\newcommand{\rdual}[1]{#1^-}
\newcommand{\ov}[1]{\overline{#1}}
\newcommand{\comp}{\circ}
\newcommand{\rowspace}{\mathcal{R}}
\newcommand{\colspace}{\mathcal{C}}
\newcommand{\mrm}[1]{\mathrm{#1}}
\newcommand{\bydef}{\stackrel{\mathrm{def}}{=}}
\newcommand{\set}[2]{\{#1\mid\,#2\}}
\newcommand{\lres}{\backslash}
\newcommand{\lresop}{\mathop{\scriptstyle\stackrel{\mathrm{op}}{\backslash}}\,}
\newcommand{\resop}{\mathop{\scriptstyle\stackrel{\mathrm{op}}{/}}\,}

\newcommand{\sh}{^{\sharp}}

\newcommand{\sil}{{\iota_{\ell}}}
\newcommand{\sir}{{\iota_r}}
\def\<#1,#2>{\langle#1\mid #2\rangle}
\newtheorem{thm}{Theorem}
\newtheorem{prop}[thm]{Proposition}
\newtheorem{lem}[thm]{Lemma}
\newtheorem{cor}[thm]{Corollary}
\theoremstyle{definition}
\newtheorem{exmp}[thm]{Example}
\newtheorem{rem}[thm]{Remark}
\date{December 20, 2002. Revised July 25, 2003}
\title{Duality and Separation Theorems in Idempotent Semimodules}
\thanks{This work was partially supported by the European Community
Framework IV program through the research network ALAPEDES (``The
Algebraic Approach to Performance Evaluation of Discrete Event Systems'').}
\author{Guy Cohen}
\address{Guy.Cohen@mail.enpc.fr : 
Cermics-Enpc, 77455 Marne-La-Vall\'{e}e, cedex 2, France.}
\author{St\'ephane Gaubert}
\address{Stephane.Gaubert@inria.fr : 
Inria-Rocquencourt, 78153 Le Chesnay cedex, France.}
\author{Jean-Pierre Quadrat}
\address{Jean-Pierre.Quadrat@inria.fr :
Inria-Rocquencourt, 78153 Le Chesnay cedex, France.}

\begin{document}
\begin{abstract}
We consider subsemimodules and convex subsets of semimodules
over semirings with an idempotent addition. We
introduce a nonlinear projection on subsemimodules:
the projection of a point is the maximal approximation
from below of the point in the subsemimodule.
We use this projection to separate 
a point from a convex set.
We also show that the projection minimizes the analogue of Hilbert's
projective metric.
We develop more generally a theory of
dual pairs for idempotent semimodules.
We obtain as a corollary duality
results between the row and column spaces of matrices
with entries in idempotent semirings.
We illustrate the results by showing
polyhedra and half-spaces over the max-plus
semiring.
\end{abstract}
\maketitle
\section{Introduction}
In this paper, we study semimodules over semirings
whose addition is idempotent, that we
call \new{idempotent semimodules}. 

A typical example of semiring
with an idempotent addition is the 
\new{max-plus semiring}, $\rmax$, which is the set
$\R\cup\{-\infty\}$, equipped with 
the addition $(a,b)\mapsto \max(a,b)$
and the multiplication $(a,b)\mapsto a+b$.
We shall also consider the \new{completed max-plus semiring},
$\rmaxb$, which is obtained by adjoining
to $\rmax$ a $+\infty$ element.
The \new{Boolean semiring} $\B$ is a
subsemiring of $\rmax$ and $\rmaxb$
(obtained by keeping only the zero element,
$-\infty$, and the unit element, $0$). 

Idempotent semimodules include a number of familiar examples.
For instance, the set of convex functions defined on a vector space
can be thought of as a semimodule over the max-plus semiring.
Another familiar class of idempotent semimodules
consists of sup-semilattices with a bottom element, which
coincide with semimodules over the Boolean semiring.

The study of idempotent analogues of linear
algebraic structures has a long history.
Early works, motivated by problems from scheduling theory, graph theory,
or dynamic programming, include \citename{cuning61} 
\citeyear{cuning61,cuning62},
\cite{H61}, \cite{yoeli61},
~\citename{V63}
 \citeyear{V63,vorobyev67,V70},
\citename{romanovski} \citeyear{romanovski},
\citename{C71a} \citeyear{C71a,C79},
~\cite{Zimmermann.K},
\citename{gondran84} \citeyear{gondran77,gondran84}.
The idempotent semimodules
that we study here were already apparent in~\cite{K65}.

More recently, the interest for idempotent semimodules arose
from the development of the max-plus algebraic approach 
to optimal control 
and asymptotic analysis~\cite{Maslov73,maslov92,maslovkolokoltsov95},
\citename{litvinov98} \citeyear{litvinov98,litvinov00},
and to discrete event systems~\cite{cohen85a,bcoq,maxplus97,ccggq99}.
See~\citename{cuning}
\citeyear{cuning,cuning95b}
\cite{Zimmermann.U,kim82,cao84,golan92,guna96,gondran02},
for more background.
Other works, dealing specially with semimodules,
are~\cite{wagneur91}, \citename{CGQ96a} \citeyear{CGQ96a,CGQ97},
\cite{nuclear}.

In this paper, we give Hahn-Banach type
theorems for complete idempotent semimodules
(the notion of completeness is defined in terms
of the natural order of the semimodule).
We show that a universal separation result holds
(Theorem~\ref{th-separ} below),
without any additional assumptions on the semimodule
or on the semiring, if one takes
as a nonlinear dual space an opposite semimodule. 
To recover a separation theorem involving a linear dual space,
we study more generally dual pairs, similar to the ones
that arise classically in the theory
of topological vector spaces: a \new{predual pair} consists of two 
complete semimodules $X,Y$,
equipped with a bilinear continuous
pairing $\<\cdot,\cdot>$, and a \new{dual pair}
is a predual pair which separates points
(see~\S\ref{sec-dual}).
We introduce a Galois connection $X\to Y, x\mapsto \ldual{x}$,
$Y\to X, y\mapsto\rdual{y}$, which yields
anti-isomorphisms between the lattices
of the elements of $X$ and $Y$ which
are closed for this correspondence.
For instance, when $X=\rmaxb^{{n\times 1}}$ is the semimodule
of $n$-dimensional column vectors over the
completed max-plus semiring $\rmaxb$,
$Y=\rmaxb^{1\times n}$, and $\<y,x>=\max_{1\leq i\leq n} (y_i + x_i)$, 
all elements of $X$ and $Y$ are closed, and
the conjugation operation is simply $\ldual{x}=(-x)^\transpose$
and $\rdual{y}=(-y)^\transpose$ where $\transpose$ denotes the transposition.
For a class of idempotent semirings that we
call reflexive, we show that dual pairs satisfy 
a more familiar, linear, geometric Hahn-Banach theorem,
which has the following form (see Theorem~\ref{t-hb} below):
if $V$ is a complete subsemimodule of $X$,
if $x\in X$ but $x\not\in V$, then there 
exist elements $y,z\in Y$ such that
\begin{align}
\<y,v>=\<z,v>,\forall v\in V,\quad\mrm{ and }\quad \<y,x>\neq \<z,x>  \enspace.
\label{e-hb}
\end{align}
\begin{figure}[hbtp]
\begin{center}
\includegraphics[scale=0.8]{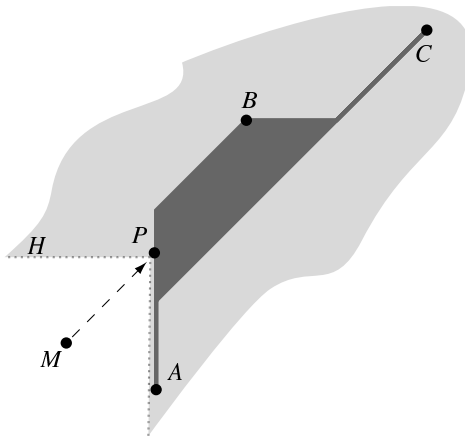}
\caption{Separation of the convex $ABC$ and the point $M$ by the half-space $H$. }
\label{fig:separe}
\end{center}
\end{figure}
The separating pair $(y,z)$ is nothing but 
the pair of conjugates $(\ldual{x}, \ldual{P_V(x)})$,
where $P_V(x)$ is the best approximation
from below of $x$ by an element of $V$.
Since $P_V(x)$ minimizes an analogue of Hilbert's projective metric,
\eqref{e-hb} is similar to the separation property in Euclidean spaces,
where $P_V$ is the orthogonal projector on $V$
and the vector $(x,P_V(x))$ gives the direction
orthogonal to a separating hyperplane.
The key discrepancy, by comparison with
vector spaces, is that one 
needs {\em pairs} of linear forms to separate a point
from a subspace, or more generally, from a convex set.
The affine form of the separation theorem is illustrated 
in Fig.~\ref{fig:separe}, which shows a max-plus
polyhedron generated by three extremal points, $A,B,C$,
a point $M$ which does not belong to the polyhedron,
together with a half space $H$ (in light gray) which
contains the polyhedron, but not the point $M$.
The half-space is obtained from the projection $P$ of $M$. 
See~\S\ref{sec-2d} for details. 

The present idempotent Hahn-Banach theorem extends several
earlier results. The first theorem of this kind seems to have
been proved by~\citename{zimmerman77}~\citeyear{zimmerman77},
for closed convex subsets of $\sK^n$, where $\sK$ is a 
semiring with an idempotent addition, satisfying
some axioms which hold when $\sK=\rmax$. 
A similar result was proved
by~\citename{shpiz} \citeyear{shpiz}
under more general assumptions
on the semiring, and in an infinite dimensional
context, but assuming that the point to separate has invertible
coordinates. The present Hahn-Banach theorem
holds under more general assumptions, 
and yields direct explicit formul\ae\ for separating hyperplanes.
This generality is possible
because we work in {\em complete} ordered structures.
In the case of the max-plus semiring, 
this means that the coefficients of the separating half-spaces
that we build can take the $+\infty$ value, so that these
half-spaces need not be closed for the usual topology.
Hence, our results apply even to some convex
subsets which are not closed,
see the example at the end of Remark~\ref{rem-rem} below.
In~\cite{singer02}, we apply the present results to convex functions
over the max-plus semiring, and recover in particular
a separation theorem \`a la Zimmerman for closed
convex sets. 
The spirit of the present work is also 
very close to that of the theory developed by
Litvinov, Maslov, and Shpiz
~\citeyear{litvinov00},
who establish idempotent analogues of several
classical theorems of functional analysis.
The representation theorem for linear
forms (Corollary~\ref{cor-2} below) and the related
analytic form the Hahn-Banach theorem (Corollary~\ref{cor-hb-ana}
below) are extensions of the corresponding
results of~\cite{litvinov00}.
Finally, we note that a preliminary version 
of the present results appeared in~\cite{CGQ00}.

We thank V. Kolokoltsov,
who suggested~\citeyear{kolokoltsov99} to the second author
the interest of revisiting max-plus residuation theory
with a Galois connection point of view: the present
work illustrates the fruitful
character of this idea, which is also applied to different problems
in~\cite{AGK00}.
We thank M.~Akian, P.~Lotito, E.~Mancinelli, I.~Singer and E.~Wagneur,
for useful discussions. We also thank the referees for their careful
reading and detailed comments. 
\section{Preliminaries}\label{sec-prelim}
\subsection{Complete Ordered Sets and Residuated Maps}
We first recall some classical notions
about ordered sets and residuated maps. 
See~\cite{birkhoff40,dubreil53,blyth72}
for more details.

By \new{ordered set}, we will mean throughout
the paper a set equipped with a partial order
relation.
For any subset $X$ of an ordered set $(S,\leq)$, we denote by $\supp X$
(resp. $\inff X$) the least upper bound (resp. greatest lower bound)
of $X$, when it exists. When $\supp X$ (resp. $\inff X$) belongs to $X$, we say
that $\supp X$ (resp. $\inff X$) is the \new{top} (resp. \new{bottom}) element
of $X$, and we write $\maxx X$ (resp. $\minn X$)
instead of $\supp X$ (resp. $\inff X$). 
We say that an ordered set $(S,\leq)$
is \new{complete} if any subset $X\subset S$
has a least upper bound. Then, $S$ has a bottom
element, $\minn S=\supp\varnothing$,
$S$ has a top element, $\maxx S=\supp S$,
and the greatest lower bound of a subset $X$ 
of $S$ is given by $\inff X=\supp\set{y\in S}{y\leq x,\; \forall x\in X}$,
so that $S$ is a complete lattice. 

If $(S,\leq)$ and $(T,\leq)$ are ordered sets, 
we say that a map $f:S\to T$ is
\new{monotone} if $s\leq s'\implies f(s)\leq f(s')$.
We say that $f$ is \new{residuated}
if there exists a map $f\sh: T\to S$ such that 
\begin{equation}
f(s) \leq t \iff s\leq f\sh(t) \enspace .
\label{e-def-res}
\end{equation}
The map $f$ is residuated
if, and only if, for all $t\in T$,
$\set{s\in S}{f(s)\leq t}$ has a top element.
Then,
\begin{align}
f\sh(t)&=\maxx\set{s\in S}{f(s)\leq t},\quad \forall t\in T \enspace,
\label{e-d-dres1}
\end{align}
which shows in particular that $f\sh$ is monotone.
If $(X,\leq)$ is an ordered set,
we denote by $(X\op,\leqop)$
the \new{opposite} ordered
set, for which $x\leqop y\iff x\geq y$.
Due to the symmetry of the defining property~\eqref{e-def-res},
it is clear that $f:S\to T$ is residuated if,
and only if, $f\sh: T\op\to S\op$
is residuated. In particular, if $f$ is residuated,
\begin{align}
f(s)&=\minn\set{t\in T}{s\leq f\sh(t)},\quad \forall s\in S \enspace,
\label{e-d-dres2}
\end{align}
and $f$ is monotone.
One also checks that $f$ is residuated
if, and only if, it is monotone,
and there exists a monotone
map $f\sh: T\to S$ such that
\begin{equation}
\label{e-res}
f\comp f\sh \leq I_T,\qquad
f\sh\comp f \geq I_S \enspace ,
\end{equation}
where $I_X$ denotes the identity map on a set
$X$. Then, $f$ and $f\sh$ satisfy~\eqref{e-d-dres1},\eqref{e-d-dres2}.

When $S,T$ are complete ordered
sets, residuated maps can be
characterized as follows.
Consider the following property,
for a map $f:S\to T$:
\begin{align}
\forall U\subset S,
\quad f(\supp U) = \supp f(U),
\quad\mrm{ where }\quad
f(U)=\set{f(x)}{x\in U}
\enspace .
\label{e-def-sci}
\end{align}
This implies in particular that $f$ is
monotone, and that $f(\minn S)=\minn T$
(take $U=\varnothing$ in~\eqref{e-def-sci}).
We shall say
that $f$ is \new{continuous} if it satisfies~\eqref{e-def-sci}.
(The term ``continuous'' can be related
to Scott topology~\cite{gierzETAL}.)
We get:
\begin{lem}
\label{lem-1}
If $(S,\leq)$ and $(T,\leq)$ are complete ordered sets,
then, a map $f:S\to T$ is residuated if, and only if,
it is continuous.\qed
\end{lem}
(See~\cite[Th.~5.2]{blyth72}, or~\cite[Th.~4.50]{bcoq}
for a proof.) 
By symmetry, if $(S,\leq)$ and $(T,\leq)$ are complete ordered sets,
and if $f$ is residuated, then, $f\sh: T\op\to S\op$, is continuous,
which means that:
\begin{equation}
f\sh(\inff U)= \inff f\sh(U) \enspace, \forall U\subset T \enspace .
\label{e-inf-morphism}
\end{equation}
We warn the reader that when $S=T=\R\cup\{\pm\infty\}$, 
a monotone map $f:S\to T$ is continuous (in the sense
of~\eqref{e-def-sci}) if, and only if, it is {\em lower semi-continuous}
in the ordinary sense and fixes $-\infty$,
whereas a monotone map $g:T\op \to S\op$ is continuous
if, and only if, it is {\em upper semi-continuous} in the ordinary sense
and fixes $+\infty$.

Using the monotonicity of $f$ and $f\sh$, together
with~\eqref{e-res}, we easily get that
\begin{subequations}
\label{e-invol}
\begin{gather}
f\comp f\sh\comp f=f  \enspace ,\label{e-invol1}\\
f\sh\comp f\comp f\sh  =f\sh \enspace ,\label{e-invol2}\\
f\sh\comp g\sh=(g\comp f)\sh  \enspace,\label{e-compose}
\end{gather}
\end{subequations}
where $g$ is a residuated map
from $T$ to some ordered set.
It is not difficult to check that
\begin{subequations}
\begin{align}
  f \text{\ is injective} \Leftrightarrow  f\sh \comp f =I_S 
 \Leftrightarrow  f\sh
\text{\ is surjective,} \label{e-injsurj}\\
  f\text{\ is surjective} \Leftrightarrow  f \comp f\sh =I_T 
 \Leftrightarrow  f\sh \text{\ is injective.}
\end{align}
\end{subequations}
Moreover, if $\{f_i\}_{i \in I}$ is an arbitrary family
of residuated maps from a complete ordered set $S$ to
a complete ordered set $T$, 
\begin{equation}
(\supp_{i\in I} f_i)\sh =\inff_{i\in I} f_i\sh \enspace ,
\label{e-supres}
\end{equation}
where the $\supp$ and $\inff$ are taken pointwise.
\subsection{Semimodules over Idempotent Semirings}
In the sequel, $(\sK,\oplus,\otimes,\zero,\unit)$
denotes a semiring whose addition is idempotent
(i.e.\/ $a\oplus a=a$), and \(\zero\) and \(\unit\) are the neutral elements for \(\oplus\) and \(\otimes\), respectively. We shall adopt the usual
conventions, and write for instance $ab$ instead
of $a\otimes b$.
An idempotent commutative monoid $(S,\oplus,\zero)$ 
can be equipped with the \new{natural} order relation,
$a\leq b\Leftrightarrow a\oplus b=b$, for which $a\oplus b=\supp\{a,b\}$
and $\zero=\minn \sK$. We say that the semiring $\sK$ is
\new{complete} if it is complete as a naturally ordered
set, and if the left and right multiplications,
$L^\sK_a, R^\sK_a: \sK\to \sK$, $L^\sK_a(x)=ax$, $R^\sK_a(x)=xa$,
are continuous.

A (right) $\sK$-\new{semimodule} $X$
is a commutative monoid $(X,\oplus,\zero)$,
equipped with a map $X\times \sK\to X$, $(x,\lambda)
\to x\lambda$ (right action), that satisfies
\begin{subequations}
\begin{gather}
x(\lambda \mu)=(x\lambda)\mu \enspace,\label{e-lefta}\\
 (x\oplus y)\lambda = x\lambda\oplus
y\lambda\;,\quad
x(\lambda \oplus \mu) = x\lambda\oplus
x\mu \enspace,\label{e-dist0}\\
x\zero = \zero\enspace,\label{e-absright}\\ 
x\unit=x \enspace,\label{e-unit}
\end{gather}
\end{subequations}
for all $x,y\in X$, $\lambda,\mu\in \sK$. 
Since $(\sK,\oplus)$ is idempotent, 
$(X,\oplus)$ is idempotent:
\[
x\oplus x = x
\]
(it follows from~\eqref{e-dist0}
and~\eqref{e-unit} that $x=x\unit=
x(\unit\oplus\unit)=x\unit\oplus x\unit=x\oplus x$).
Axiom~\eqref{e-absright}
may be we rewritten more explicitly as $x\zero_\sK=\zero_X$.
It implies that 
\begin{align}
\zero_X \lambda = \zero_X \enspace .
\label{e-dualabs}
\end{align}
Indeed, for any $x\in X$, $\zero_X \lambda=(x\zero_\sK)
\lambda=x(\zero_{\sK}\lambda)=x\zero_{\sK}=\zero_{X}$,
using~\eqref{e-lefta} and the fact that
$\zero_{\sK}$ is absorbing for the product of $\sK$.

The notion of left $\sK$-semimodule
is defined dually.  Throughout the paper,
all the semimodules that we shall consider
will be over idempotent semirings.
We shall also consider $\sK$-\new{bisemimodules}:
a bisemimodule is a set equipped
with two, right and left, $\sK$-semimodule structures, such that
the right and left actions commute. In particular,
an idempotent semiring $\sK$ is a $\sK$-bisemimodule
if one take as left and right actions the semiring product
$(a,b)\mapsto ab$. 

When $\sK$ is a complete idempotent semiring,
we say that a right $\sK$-semimodule $X$ is \new{complete}
if it is complete as a naturally ordered set,
and if, for all $v\in X$ and
$\lambda\in \sK$,
the left and right multiplications,
$R^X_{\lambda}:\;X\to X$, $x\mapsto x\lambda$
and $L^X_{v}:\;\sK\to X$, $\mu\mapsto v\mu$,
are both continuous. 
Complete left $\sK$-semimodules
and complete $\sK$-bisemimodules are defined
in a similar way.  In the sequel, all semimodules
will be right semimodules, unless otherwise
specified. We shall also use the notion
of \new{linear map} (as usual, a map between semimodules
is linear if it preserves finite sums and commutes with the action).
\begin{exmp}[Free Complete Semimodules and Semimodules of Functions]\label{ex-free}
Let $\sK$ denote a complete idempotent semiring.
A {\em free complete right $\sK$-semimodule}
is of the form $\sK^I$ for some arbitrary set $I$:
the elements of $\sK^I$ are functions $I\to \sK$,
and $\sK^I$ 
is equipped with the addition
$(a,b)\mapsto a\oplus b, \; (a\oplus b)(i)=a(i)\oplus b(i)$,
and the action $(a,\lambda)\mapsto a\lambda, \;
(a\lambda)(i)=a(i)\lambda$, for
all $a,b\in \sK^I, \lambda\in \sK$.
By considering the action 
$(a,\lambda)\mapsto \lambda a,(\lambda a)(i)
=\lambda a(i)$, one can see $\sK^I$
as a left semimodule.

The semimodule \(\rmaxb^{n\times 1}\), evoked in the introduction,  is an example of a free complete right \(\rmaxb\)-semimodule.
\label{ex:funct}
 Another example in the same category, to which we will return from time to time in this paper, is the set  $\rmaxb^{\sU}$ of functions from a set \(\sU\) to \(\rmaxb\), with the pointwise supremum as \(\oplus\) operation and the conventional addition of a real constant as (left or right) action. This semimodule (that we refer to as \(\sF\) for short in the sequel) is complete.
\end{exmp}
In a complete semimodule $X$, we define,
for all $x,y\in X$ and $\lambda\in \sK$,
\begin{subequations}
\begin{align}
  x\lres y &\bydef (L_x^X)\sh(y) = \maxx\set{\lambda\in \sK}{x\lambda \leq y} \enspace ,\label{e-moins}\\
  x/ \lambda &\bydef (R_{\lambda}^X)\sh(x)=\maxx
\set{y\in X}{y\lambda\leq x}  
\label{e-zero}
\end{align}
\end{subequations}
(recall our convention to write $\maxx$, instead of $\supp$,
to emphasize the fact the the supremum belongs to the set).
Paraphrasing the definition of residuated maps, 
\begin{align}
\label{e-1}
x\lambda \leq y\iff \lambda \leq x\lres y 
\iff x\leq y/\lambda \enspace .
\end{align}
The residuation formul\ae~\eqref{e-res},
~\eqref{e-inf-morphism}, ~\eqref{e-invol}
and~\eqref{e-supres} yield
\begin{subequations}
\label{e-globalres}
\begin{alignat}{2}
&x(x\lres y)  \leq y \;, &&\qquad(x/\lambda)\lambda 
\leq x \;, \label{e-0} \\
&(x\lres y)\lambda \leq  x\lres(y \lambda ) \;,&&\qquad x(\lambda / \mu)
  \leq (x \lambda)/\mu\;,\label{e-lnum}\\
&x\lres (x\lambda)  \geq  \lambda \;, &&\qquad (x\lambda)/\lambda \geq x \;, \\
& x\lres (\inff U) =  \inff (x\lres U)
\;, &&\qquad (\inff  U) / \lambda  =  \inff (U / \lambda)\;,\label{e-inf}\\
& x(x\lres (x\lambda)) =  x\lambda \;, &&\qquad((x\lambda)/\lambda)\lambda 
  = x\lambda \;, \label{e-xxxy}\\
&  x\lres (x(x\lres y))  =  x \lres y 
\;, && \qquad((x/\lambda)\lambda)/\lambda = x / \lambda\;, \\
&\lambda\lres (x\lres z) =(x\lambda)\lres z\;, &&\qquad(x/\mu)/\lambda  
 = x/(\lambda\mu)\label{e-comp}\;,\\
& (\supp U)\lres y  =\inff (U\lres y) \;, &&\qquad x/(\supp \Lambda) 
 = \inff (x/\Lambda) \;,\label{e-11}
\end{alignat}
\end{subequations}
for all $x,y,z\in X$, $U\subset X$,
$\Lambda\subset \sK$,
where $(U\lres y)=\set{u\lres y}{u\in U}$,
$x/\Lambda=\set{x/\lambda}{\lambda\in \Lambda}$, etc.
(In the right formula~\eqref{e-lnum} and in the left formula~\eqref{e-comp},
$\sK$ is seen respectively as a right and left
semimodule over itself.)
Finally, if $X$ is a bisemimodule and $\mu,\nu\in \sK$,
the maps $x\mapsto x\lambda$
and $x\mapsto \nu x$ commute, hence, by~\eqref{e-compose},
their residuated maps commute, which means that
\begin{align}
(\nu \lres x)/\mu &=  \nu \lres (x /\mu) \enspace.
\label{e-cm}
\end{align}
Since there is no ambiguity, we may simply write $\nu
\lres x/\mu$ for~\eqref{e-cm}.

\begin{rem}
Note that \eqref{e-moins} is dual of the definition~(5.1)
in \citename{litvinov00}; the latter requires the assumption that
the action of  vectors on scalars  satisfies
$x(\inff_{\lambda \in \Lambda}\lambda ) =\inff_{\lambda \in \Lambda} (x\lambda)$ --- see \citename[Eq.~(4.7)]{litvinov00}
which is written for right action of vectors on scalars ---
whereas, in this paper, we stick to the more natural assumption that this property holds with \(\supp\) instead of \(\inff\): 
this is the case for instance if the underlying semiring
is a semiring  of formal series, or of matrices, over a complete
idempotent semiring.
\end{rem}
\subsection{Opposite Semimodules}
If $X$ is a complete right $\sK$-semimodule, 
we call \new{opposite} semimodule of $X$ the {\em left} 
$\sK$-semimodule
$X\op$ with underlying set $X$, addition $(x,y)\mapsto
\inff\{x,y\}$ (the $\inff$ is for the natural order of $X$)
and left action $\sK\times X\to X$, $(\lambda, x)\to 
x/\lambda$.
For clarity, we shall sometimes
denote by $(\lambda,x)\mapsto \lambda \dotop x=x/\lambda$
the left action of $X\op$.
That $X\op$ is a complete semimodule
follows from formul\ae~\eqref{e-inf},~\eqref{e-comp}, and~\eqref{e-11}. 
In particular,~\eqref{e-comp} yields
\begin{align}
(\lambda \mu)\dotop x&= x/(\lambda \mu)
=(x/\mu)/\lambda
= \lambda\dotop (\mu\dotop x)
\enspace, \label{e-lefta2}
\end{align}
for all $\lambda,\mu\in \sK$
and $x\in X\op$,
which shows why $X\op$ must be considered
as a {\em left} rather than a {\em right} semimodule.
Indeed, considering $(x,\lambda)\mapsto x/\lambda$ as
a right action would require
the property symmetrical to~\eqref{e-lefta2}
to hold, that is, by~\eqref{e-lefta}, $x/(\lambda\mu)=(x/\lambda)/\mu$,
but this property need not hold for a semimodule $X$
over a {\em noncommutative} semiring $\sK$.

Denoting by $\lresop$ and $\resop$
the residuated operations built from $\dotop$,
we get from~\eqref{e-1}, 
\begin{subequations}\label{e-op}
\begin{align}
\lambda \lresop x= (L_\lambda^{X\op})\sh(x)&= \minn\set{y\in X}{y/\lambda \geq x}  
= x\lambda \enspace, \label{e-nd1}\\
x\resop y = (R_y^{X\op})\sh(x)&= \maxx\set{\lambda\in \sK}{y/\lambda \geq x} 
= x\lres y \enspace .\label{e-nd2}
\end{align}
\end{subequations}
Eqn~\eqref{e-nd1} is an involutivity property: the residuated
law of the residuated law of the right action of $X$ is the
right action of $X$ itself. Therefore,
\begin{prop}
For all complete $\sK$-semimodules $X$, $(X\op)\op=X$.
\qed 
\end{prop}
\section{Nonlinear Projectors, Universal Separation Theorem
and Hilbert projective metric}
\label{sec-proj}
\subsection{Nonlinear Projector}
\label{sec-separ}
Let $V$ denote a \new{complete subsemimodule} of a complete semimodule $X$
over a complete idempotent semiring $\sK$, i.e., 
a subset of $X$ that is stable by arbitrary sups
and by the action of scalars.
We call \new{canonical projector} on $V$ the map
\[
P_V: X\to X,\quad P_V(x) = \maxx\set{v\in V}{v\leq x}	\]
(the least upper bound of $\set{v\in V}{v\leq x}$ 
belongs to this set by definition of complete subsemimodules).
It is readily seen that $P_V^2=P_V$ and that $P_V(X)=V$.
We say that $W$ is a \new{generating family}
of a complete subsemimodule $V$ if any
element $v\in V$ can be written as 
$v=\supp\set{w \lambda_w }{w\in W}$,
for some $\lambda_w\in \sK$.
\begin{thm}[Projector Formula]
If $V$ is a complete subsemimodule of $X$ with generating
family $W$, then
\begin{align}
P_V(x) = \supp_{w\in W} w (w\lres x) \enspace .
\label{e-proj}
\end{align}
\end{thm}
\begin{proof}
We can write $P_V(x)=\supp_{w\in W}w\lambda_w$,
for some $\lambda_w\in \sK$. From $P_V(x)\leq x$,
we get $w\lambda_w\leq x$, 
or, equivalently, $\lambda_w\leq w\lres x$.
This shows that $P_V(x)\leq \supp_{w\in W} w(w\lres x)$.
But, $\supp_{w\in W} w(w\lres x)$
is an element of $V$, which, by~\eqref{e-0}, 
is less than or equal to $x$. This proves~\eqref{e-proj}.
\end{proof}
We may rewrite~\eqref{e-proj} as 
\(
P_V = \supp_{w\in W} P_w \enspace, 
\)
where $P_w$ denotes the projector on the ``one dimensional''
space $w\sK$. Similar formul\ae\ for
the projector appeared in~\cite{moller88}.
\begin{prop}[Dual characterization of the projector]
Let $V\subset X$ denote a complete subsemimodule
with generating family $W$. Then,
\begin{equation}
  P_V(x) =  \minn\set{z\in X}{ w\lres z \geq w\lres x,\;\forall w\in W} \label{e-2}\;.
\end{equation}
\end{prop}
\begin{proof}
Since $w\lres z\geq w\lres x\iff z\geq w(w\lres x)$,
this follows from~\eqref{e-proj}.
\end{proof}
\begin{exmp}\label{ex:conv}
We return to the \(\rmaxb\)-semimodule \(\sF\) introduced at Example~\ref{ex:funct} and discuss the application of previous results in this section. First of all, observe that 
\begin{displaymath}
 \forall f, g\in \sF,\enspace f\lres g = \inf_{u\in\sU}\big(g(u)-f(u)\big)\;,
\end{displaymath}
with the convention here that \(+\infty - \infty=+\infty\),
since in any complete idempotent 
semiring~\(\sK\), \(\minn \sK\lres \minn \sK=\maxx \sK\lres \maxx \sK = \maxx \sK\).
(Observe however that \((\minn \sK)(\maxx \sK) = \zero (\maxx \sK)=\zero= \minn \sK\), which translates, in \(\rmaxb\), as \(-\infty+\infty=-\infty\), so that this ``rule'' written in conventional notation is ambiguous, and one must keep in mind what are the correct algebraic operations hidden behind the conventional notation to apply the rule correctly.)

Assume now that $\sU$ is a locally convex
topological vector space and consider the complete subsemimodule \(V\) generated by the set \(W\) of continuous linear functions over~\(\sU\).
This semimodule consists of the identically $-\infty$ function
over $\sU$, and of the l.s.c.~convex functions over \(\sU\)
which do not take the value $-\infty$.
For any \(f\in\sF\), \(P_V(f)\), as defined in \S\ref{sec-separ}, is the classical l.s.c.~convex hull of \(f\). For \(w\in W\), 
\begin{displaymath}
w\lres f=\inf_{u\in\sU}\big(f(u)-w(u)\big)=-\sup_{u\in\sU}\big(w(u)-f(u)\big)\;,
\end{displaymath}
which coincides, up to a change of sign, with the Legendre-Fenchel transform \(f\fench\) of \(f\) evaluated at \(w\). Eqn~\eqref{e-proj} then yields
\begin{displaymath}
P_V(f)(\cdot)=\supp_{w\in W}w(w\lres f)=\supp_{w\in W}\big( w(\cdot)-f\fench(w)\big)\;,
\end{displaymath}
that is to say, the l.s.c.~convex hull of \(f\) is the Legendre-Fenchel transform of the Legendre-Fenchel transform of \(f\).

Finally, Eqn~\eqref{e-2} says that the l.s.c.~convex hull of \(f\) is the least function~\(g\) in \(\sF\) such that \(g\fench\) is less than, or equal to, \(f\fench\) (pointwise).
\end{exmp}
\subsection{Universal Separation Theorem}
\begin{thm}[Universal Separation Theorem]
\label{th-separ}
Let $V\subset X$ denote a complete subsemimodule, 
and let $x\in X$. Then, 
\begin{subequations}
\label{e-separ}
\begin{align}
\forall v\in V,&\quad v\lres P_V(x) = v\lres x \; ,
\label{sa}\\\intertext{and}
x\in V&\iff x\lres P_V(x) = x\lres x\; .
\label{sb}
\end{align}
\end{subequations}
\end{thm}
Seeing $y\lres x$ as a ``scalar product'',
Eqn~\eqref{sa} says that the vector
$(x,P_V(x))$ is ``orthogonal'' to the semimodule $V$,
and~\eqref{sb} shows that the 
``hyperplane'' $\set{y}{y\lres P_V(x)=y\lres x}$
separates $x$ from $V$, if and only if $x\not\in V$.
This terminology will be justified in~\S\ref{sec-dual}.
\begin{proof}
Since, by definition, the $\minn$ in~\eqref{e-2}
belongs to the set, we have that $v\lres P_V(x)\geq v\lres x$, for all $v\in V$.
Using $P_V(x)\leq x$ and the monotonicity of $y\mapsto
v\lres y$, we get
the reverse inequality, which shows~\eqref{sa}.
If $x\in V$, then $P_V(x)=x$, and $x\lres  P_V(x) = x\lres x$,
trivially. 
Conversely, 
if $x\lres  P_V(x) = x\lres x$, we have, by~\eqref{e-1},
that $P_V(x) \geq x(x\lres x)$, and, by~\eqref{e-xxxy},
that $x(x\lres x)=x$,
which shows that $P_V(x)\geq x$. Since $P_V(x)\leq x$,
we have $x=P_V(x)\in V$.
\end{proof}
\begin{rem}
The separating set $H=\set{v\in V}{v\lres P_V(x) = v\lres x}$
is a semimodule. Indeed, by (\ref{e-11}) it is stable by addition
 and  (\ref{e-comp}) shows that it is stable by scalar action.
\end{rem}
\begin{rem}
According to the previous remark, it is sufficient to check \eqref{sa} only for \(v\) ranging in a generating subset \(W\) of \(V\).
\end{rem}
\begin{exmp}
For the semimodule \(\sF\) introduced at Example~\ref{ex:funct} and the subsemimodule~\(V\) of l.s.c.~convex functions generated by the subset \(W\) of continuous linear functions as discussed at Example~\ref{ex:conv}, the equality~\eqref{sa} (restricted to \(v\in W\) as observed in the previous remark) of the Separation Theorem says that the Legendre-Fenchel transform of any function~\(f\) coincides with the Legendre-Fenchel transform of its l.s.c.~convex hull. As for \eqref{sb}, observe first that \(f\lres f=0\) unless  \(f\) assumes only \(\pm\infty\) values (in this latter case, \(f\lres f=+\infty\)). Let us put aside this singular situation first. Then \eqref{sb} says that \(f\) coincides with its l.s.c.~convex hull at all points if and only if it is itself l.s.c.~convex. 

In the singular case, and according to \eqref{sb}, \(f\) is l.s.c.~convex if and only if \(f\lres P_V(f)=\inf_{u\in\sU}\big(P_V(f)(u)-f(u)\big)=+\infty\), that is, \(P_V(f)(u)-f(u)=+\infty\) for all \(u\). According to the rule \(-\infty+\infty=+\infty\) which applies here, this shows that \(f(u)=+\infty\) implies that \(P_V(f)(u)=+\infty\). On the other hand, if \(f(u)=-\infty\), then \(P_V(f)(u)=-\infty\) because \(P_V(f)\leq f\) pointwise. Finally, in all cases, we have reached the conclusion that \eqref{sb} says that \(f\) coincides with its l.s.c.~convex hull at all points if and only if it is itself l.s.c.~convex.
\end{exmp}
The ``scalar product''
$y\lres x$ separates points, in the
following sense:
\begin{prop}[Separation of Points]\label{prop-separpoint}
If $X$ is a complete $\sK$-semimodule,
then, for all $x,y\in X$, 
\begin{align}
(\forall z\in X,\quad x\lres z = y\lres z) \implies x = y 
\enspace .
\label{e-separpoints}
\end{align}
\end{prop}
\begin{proof}
If $x\lres z = y\lres z$ for all $z\in X$, taking
$z=x$, we get that $\unit\leq x\lres x= y\lres x$,
hence $y\leq x$. By symmetry, $x\leq y$.
\end{proof}
Finally, we note that all the above results have
dual versions for the semimodule $X\op$: they
are derived readily from~\eqref{e-op}.
For instance, if $V\subset X\op$ is a complete
subsemimodule, we define 
\begin{align}
P\op_{V}(x) &= \suppop\set{v\in V}{v\leqop x}
=  \inff\set{v\in V}{v\geq x} \enspace,
\label{e-def-pop}
\end{align}
where $\suppop=\inff$ denotes the least
upper bound associated with $\leqop$,
and the dual version of Theorem~\ref{th-separ}
reads:
\begin{thm}[Dual Separation Theorem]
Let $V\subset X\op$ denote a complete subsemimodule, 
and let $x\in X$. Then, 
\begin{subequations}
\begin{align}
\forall v\in V,&\quad P\op_V(x)\lres v = x \lres v  \; ,
\label{dualsa}\\\intertext{and}
x\in V&\iff P\op_V(x)\lres x = x\lres  x\enspace .
\label{dualsb}
\end{align}
\end{subequations}
\end{thm}
In the same way, dualizing~\eqref{e-separpoints},
we get the following separation property for points:
\begin{align}
(\forall z\in X,\quad z \lres x = z\lres y) \implies x = y 
\enspace .
\label{e-dualseparpoints}
\end{align}
\begin{rem}
It is natural to ask whether the projector
\[
Q_V(x)=\inff\set{v\in V}{v\geq x}
\]
can be defined when $V$ is a subsemimodule
of $X$, rather than a semimodule of $X\op$ as in~\eqref{e-def-pop}.
The difficulty is that $Q_V(x)$ need
not belong to $V$.
For instance, when $V\subset\rmaxb^3$
is the subsemimodule generated by the columns
of the matrix
\[
\begin{pmatrix}
0& -1\\
-1&0\\
0 & 0
\end{pmatrix} \enspace ,
\]
\[
Q_V\begin{pmatrix}-1 \\ -1 \\ 0
\end{pmatrix}
= \begin{pmatrix}-1 \\ -1 \\ 0
\end{pmatrix}
\]
does not belong to $V$.
However, in the special case
when $V$ is a complete
subsemimodule of $X$ stable by arbitrary infs,
we have $Q_V(x)\in V$, for all $x\in X$, and $Q_V$
preserves arbitrary sups, whereas
$P_V$ need not have this property.
\end{rem}
We now derive from Theorem~\ref{th-separ} a Hahn-Banach
theorem for complete convex subsets, in the spirit of~\cite{CGQ00}.
We say that a subset $C$ of a complete semimodule
over a complete semifield
$\sK$ is \new{convex} (resp.\ \new{complete convex})
if for all finite (resp.\ arbitrary)
families $\{x_i\}_{i\in I}\subset C$ and $\{\alpha_i\}_{i\in I}
\subset \sK$, such that $\supp_{i\in I} \alpha_i =\unit$,
we have that $\supp_{i\in I} \alpha_i x_i \in C$. 
Theorem~\ref{th-separ} has an immediate extension
to convex sets.
\begin{cor}[Separating a Point from a Convex Set]\label{cor-convex}
If $C$ is a complete convex subset of a complete $\sK$-semimodule $X$,
and if $x\in X$ is not in $C$, then we have
\begin{subequations}\label{e-sconvex}
\begin{align}
v\lres x \wedge e&=v\lres y \wedge \nu\; ,\qquad\forall v\in C\;,\\
x\lres x\wedge e &> x\lres y \wedge\nu   \;,
\label{e-hbc}
\end{align}
\end{subequations}
with
\begin{align}
\nu &=\supp_{v\in C}(v\lres x\wedge e)\;,
\qquad y = \supp_{v\in C}v(v\lres x\wedge e)   \;.
\label{e-formul}
\end{align}
\end{cor}
\begin{proof}
Consider the complete $\sK$-semimodule $Y=X\times \sK$
and the complete subsemimodule $V$ generated by the vectors
$(v\lambda,\lambda)$, where $v\in C$ and $\lambda \in \sK$.
It is easy to see that $(v,e)$ belongs to $V$ iff 
$v$ belongs to the complete convex set 
generated by $C$, which coincides with $C$.
When $x\not\in C$, then $(x,\unit)\not\in  V$,
and applying Theorem~\ref{th-separ}, we have that
$$(v,e)\lres (x,\unit)=(v,e)\lres P_V((x,\unit)),\quad
\forall v \in C \; .$$
$$(x,\unit)\lres (x, \unit)>(x,\unit)\lres P_V((x,\unit))\; .$$
By using this result with  
\begin{align*}(y,\nu)&=P_V((x,e))\\
&=\supp_{v\in C}(v,e)\big((v,e)\lres(x,e)\big)\\
\intertext{(thanks to \eqref{e-proj})}
&=\supp_{v\in C}(v,e)(v\lres x\wedge e)
\end{align*}
(since  $(a,\lambda)\lres (b, \lambda')= a \lres b \wedge \lambda
\lres \lambda'$),  the proof is completed.
\end{proof}
\begin{rem}\label{rem-rem}
Observe that if $x\in C$, then $P_V((x,e))=(x,e)$, $\nu=e$ and $y=x$. Moreover, if $\nu$ is invertible, then it is easy to see that $y\nu^{-1}$ belongs to $C$ and can thus be considered as the projection of $x$ onto the convex subset $C$. Indeed, setting $P_C(x)=y\nu^{-1}$ (whenever this expression is defined), the image of $C$ by $P_C$ is $C$ and $P_C\comp P_C=P_C$. 

When $\nu$ is not invertible (in $\rmaxb$, this means that $\nu=\zero$ since $\nu$ is not greater than $\unit$), we still do have a separating equation but its interpretation in terms of projection onto $C$ is missing. This happens in the following example: $X=\sK=\rmaxb$ and $C=(-\infty,+\infty]$. This $C$ is complete convex but not closed in the usual topology. Nevertheless, the previous theory still applies and we can separate $x=-\infty$ from $C$.  Calculations show that $y=\nu=-\infty$ and relations~\eqref{e-sconvex} can be checked to be true.
\end{rem}
\subsection{Generalized Hilbert projective metric}
Consider $d_H:X\times X \to \sK$ defined by $d_H(x,y)=(x\lres y)
(y\lres x)$. Observe that $d_H(x,y)=d_H(y,x)$ when $\sK$ is commutative.
When $X=\rmaxb^n$, $d_H$ is nothing but an additive
version of Hilbert projective metric,
which is the map
\[\delta_H(x,y)=\max_{1\leq i,j\leq n} \log\Big( \frac{x_i}{y_i}\frac{y_j}{x_j} \Big)\]
for $x,y$ ranging in the open positive cone of $\R^n$. When $x,y\in \R^n$, 
\[d_H(x,y)=\min_{1\leq i,j\leq n}  (x_i-y_i+y_j-x_j)=-\delta_H(\exp{x},\exp{y})\;,\]
where $\exp$ operates coordinatewise.

\begin{thm}
The map $d_H$ satisfies the following properties:
\begin{itemize}
\item[-] anti-triangular inequality (when $\sK$ is commutative):
$$d_H(x,z)\geq d_H(x,y) d_H(y,z)\; ;$$
\item[-] definiteness: 
$$d_H(x,y)=e \Rightarrow x= y   \lambda, \quad\lambda \in \sK\; ,$$ 
\item[-] nonpositiveness:
$$d_H(x,y)\leq x\lres x\enspace \text{and\ } d_H(x,y)\leq (x\lres x)\wedge( y\lres y) \text{\ when $\sK$ is commutative.} $$
\end{itemize}
\end{thm}
\begin{proof}\ \ 
\begin{itemize}
\item[-] Anti-triangular inequality: 
\begin{align*}(x\lres y)(y\lres x)(y\lres z)  
(z\lres y)&=(x\lres y) (y\lres z) (z\lres y)  
(y\lres x)
\leq (x \lres z) (z\lres x)\end{align*} 
by (\ref{e-lnum}) and (\ref{e-0}).
\item[-] Definiteness:
if $d_H(x,y)=e$ we have that
$$x=x (x\lres y)(y \lres x)\leq y (y\lres x)\leq x\;, $$hence $x=y(y\lres x)$.
\item[-] Nonpositiveness:
$$(x \lres y)  (y \lres x)\leq x\lres (y   (y\lres x))\leq x \lres x\;.$$
\end{itemize}
\end{proof}
In conventional Euclidean spaces, the projection of a point
onto a subspace minimizes the distance from that point to any point
of the subspace.
We show here that $d_H$ is maximized
by projection.

\begin{thm}\label{th:hilbert}
For all $x\in X$ and $v \in V$, where $V$ is a complete subsemimodule
of a semimodule $X$, we have that $d_H(x,v)\leq d_H(x,P_V(x))$.
\end{thm}
\begin{proof}
\begin{alignat*}{3}
d_H(x,P_V(x)) & = (x\lres P_V(x))(P_V(x)\lres x)&& && \\
      & \geq x\lres P_V(x)&& &&\qquad  \text{(because } P_V(x) \leq x) \\
      &= x\lres\Big(\supp_{v\in V}( v(v \lres x))\Big)&& &&\\
      &\geq x \lres( v(v\lres x))\;,&& \forall v\in V && \\
      & \geq (x\lres v)(v\lres x)\;,&& \forall v\in V && \qquad(\text{by (\ref{e-lnum}))} \\
      &=d_H(x,v)\;,&&\forall v\in V\; . && 
\end{alignat*} 
\end{proof}
\begin{exmp}
Once again we return to our favorite illustration described at Examples~\ref{ex:funct} and \ref{ex:conv}. For two functions \(f\) and \(g\) in \(\sF\), we consider 
\begin{displaymath}
-d_H(f,g)=\sup_{u\in\sU}\big(f(u)-g(u)\big)+\sup_{v\in\sU}\big(g(v)-f(v)\big) \enspace .
\end{displaymath}
This is a ``form factor'', which measures ``how far'' is $f-g$ from
a constant map. Indeed, when $f,g$ are finite, $-d_H(f,g)$
is nothing but the difference between the sup and the inf of of $f-g$.
Then, Theorem~\ref{th:hilbert} says that the l.s.c.~convex hull of \(f\) is, among all l.s.c.~convex functions, one which minimizes this form factor difference with \(f\) (but of course not the only one).
\end{exmp}
\subsection{A Two Dimensional Example} \label{sec-2d}
\begin{figure}[hbtp]
\begin{center}
\includegraphics[scale=0.8]{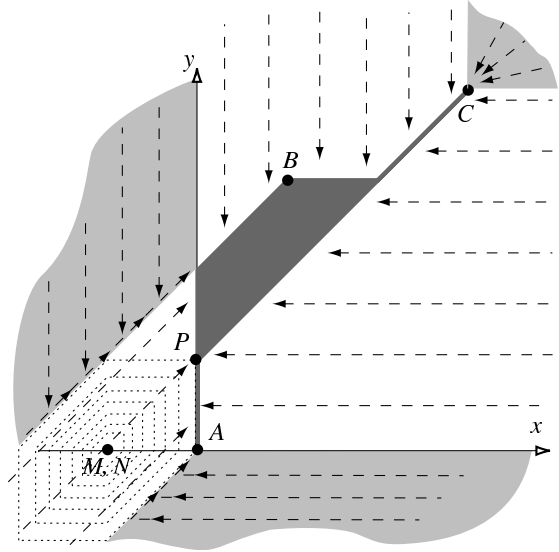}
\caption{The view in the \((x,y)\)-plane}
\label{fig:haut}
\end{center}
\end{figure}

\begin{figure}[hbtp]
\begin{center}
\includegraphics[scale=0.8]{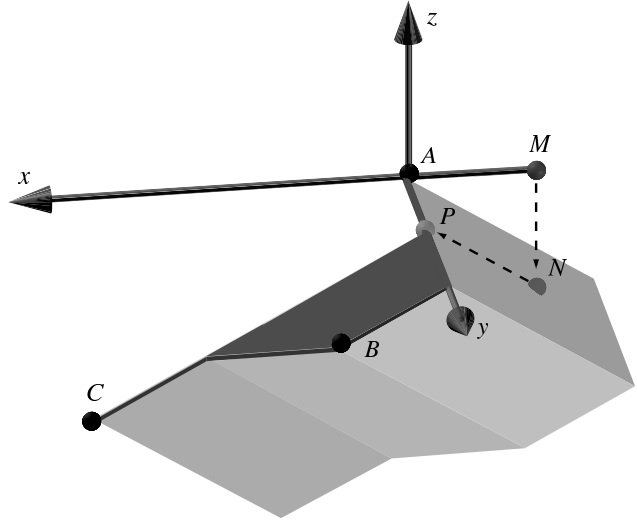}
\caption{The 3D view}
\label{fig:3D}
\end{center}
\end{figure}

\begin{figure}[hbtp]
\begin{center}
\includegraphics[scale=0.8]{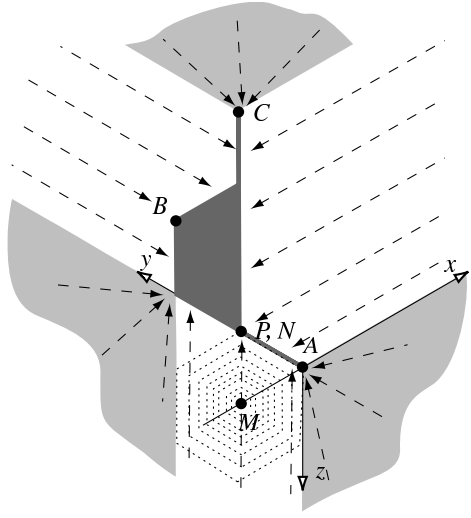}
\caption{The view of an observer located along the vector \((1,1,1)\)}
\label{fig:maire}
\end{center}
\end{figure}

We consider the convex set generated by points \(A,B,C\) of coordinates \((0,0), (1,3)\) and \((3,4)\) in \(\rmax^{2}\). Figure~\ref{fig:haut}
represents these 3~points in this space and the convex set is depicted in dark grey (notice it has two ``antennas'' ending in \(A\) and \(C\) in addition to the polygon with nonempty interior). Figure~\ref{fig:3D} is a
representation in the 3D space where (a fragment of) the subsemimodule~$V$
--- introduced in the proof of Corollary~\ref{cor-convex} ---  generated by points~\(A,B,C\) (now with coordinates \((0,0,0),(1,3,0)\) and \((3,4,0)\)) is represented. The intersection of this subsemimodule with the \((x,y)\)-plane is the convex set represented in Figure~\ref{fig:haut}. The ``cylinder'' is parallel to the vector~\((1,1,1)\). Figure~\ref{fig:maire}
is a representation of what can be seen by an observer located at a remote point along the vector \((1,1,1)\).

We now consider projecting the point \(M\) of coordinates \((-1,0)\) (in \(\rmax^{2}\)) onto the convex set. According to Remark~\ref{rem-rem}, this point is first projected on the subsemimodule~$V$ at point \(N\) of coordinates 
\((-1,0,-1)\) in \((\rmax^{3})\): indeed, this is the ``best approximation from below'' of \(M\) by an element of the subsemimodule. The reader can check this claim by using the provided explicit formul\ae~\eqref{e-proj}. Then, \(N\) is brought back to \(\rmax^{2}\) by ``normalization'' of the \(z\)-coordinate to \(0\), yielding the point~\(P\) of coordinates \((0,1,0)\). Points \(M,N,P\) are shown in the three figures. 

Relations~\eqref{e-sconvex} yield the following
\begin{align*}
\min(-1-x,-y,0)&=\min(-1-x,-y,-1), \enspace \forall (x,y) \text{\ in the convex set;} \\ \min(-1-(-1),0,0) &>\min(-1-(-1),0,-1)\enspace\text{when applied to $M$}\;.
\end{align*}
 The former equation simplifies into $\min(-1-x,-y,0)\leq -1$ which says that $-1-x\leq -1$ or $-y \leq -1$: this is the union of two half planes,
corresponding to the light grey region in Figure~\ref{fig:separe}.

Observe that in Figure~\ref{fig:haut}, points \(M\) and \(N\) are located at the same place because it turns out that they are located on the same vertical line of \(\rmax^{3}\), whereas in Figure~\ref{fig:maire}, points \(N\) and \(P\) are located at the same place: this is a general fact because normalization always implies a move in the direction in which the observer of this figure is located.

In Figures~\ref{fig:haut} and \ref{fig:maire}, several zones around
 the convex set are also shown: \begin{itemize}
\item in light grey conic zones, it turns out
 that all points project onto a particular ``extreme'' point of the
 convex set;\begin{itemize}
\item in the grey zone attached to point~\(C\) (and in the whole 
positive orthant \((x\geq 0, y\geq 0)\) as well), there is a single
 move in the \((x,y)\)-plane, that is, the projection onto the  subsemimodule 
coincides with that onto the convex set;\item in the other two grey zones,
 there are actually two moves: one caused by the projection onto the
 subsemimodule, the other one caused by normalization; this is
 materialized by dotted line arrows in Figure~\ref{fig:haut}; 
in Figure~\ref{fig:maire}, the latter move (caused by normalization) 
is not visible for reasons already explained hereabove.
\end{itemize}\item in the white 
zones of the positive orthant, as already mentioned, the moves are 
always one-phase (i.e.~horizontal); in the white zone which \(M\)
 belongs to,  the former move is vertical (thus it cannot be visualized on
 Figure~\ref{fig:haut}) 
and the latter one (normalization) is (as everywhere) along the first diagonal.
\end{itemize}

Finally, level sets of the generalized Hilbert metric are shown around point~\(M\) in those figures.
\begin{exmp}
It is useful to understand  the geometry of
affine max-plus hyperplanes of $\rmax^2$,
that we shall call {\em lines}.
The general line is defined by 
an equation of the form $$ax \oplus by \oplus c
= a'x\oplus b'y \oplus c'\;,
$$for some $a,b,c,a',b',c'\in \rmax$,
but not all the coefficients are
needed. For instance, the lines
with equations $2 x \oplus y = 1x  \oplus y \oplus 3$
and $2x \oplus  y = y \oplus 3$, coincide.
More generally, it is not difficult to see that there are $12$
generic shapes of lines, as shown
in Figure~\ref{fig:lines}. 
\begin{figure}[hbtp]
\begin{center}
\includegraphics[scale=0.8]{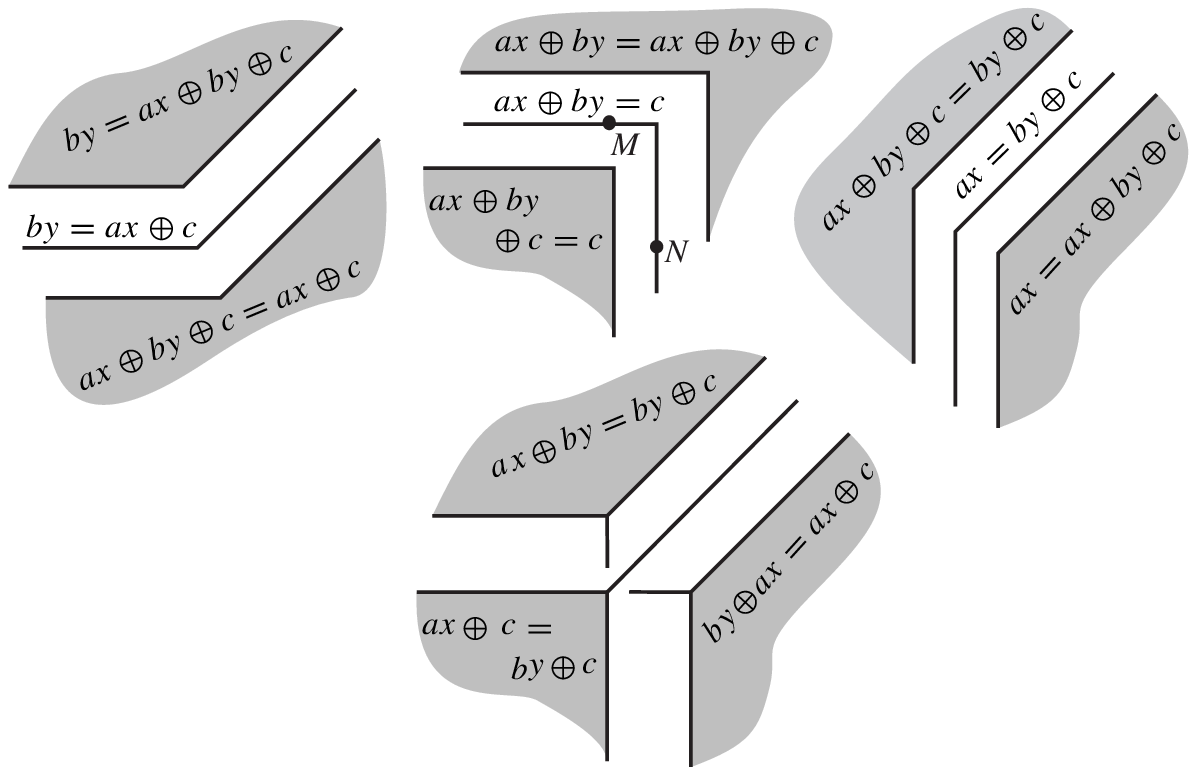}
\caption{The twelve generic lines of \(\rmax^{2}\)}
\label{fig:lines}
\end{center}
\end{figure}
Indeed, a generic line can be defined by three real numbers 
\(a,b,c\) plus a ``sign'' information, 
which tells the side of the equation in which the corresponding 
coefficients is dominant (say ``$\oplus$''
for the left hand side, ``$\ominus$'' for the right
and side, and a dot when coefficients on both
sides are equal). For instance, 
the line with equation 
\(ax\oplus c= b y\oplus c\) can be denoted \(L(\oplus a,\ominus 
b, \dot{c})\). 
This notation can be justified by introducing
the {\em symmetrized} max-plus semiring~\cite{maxplus90b,gaubert92a,bcoq}.
It is fundamental to note that a line with a dotted
coefficient has dimension $2$ in the usual
sense. Note also that half-planes are special lines,
since for instance an inequality of the form $x\geq y$ can be written as an equation $x=x\oplus y$.
\end{exmp} 
\section{Dual Semimodules and Hahn-Banach Theorems}\label{sec-dual}
\subsection{Dual and Predual Pairs}
Given a complete idempotent semiring $\sK$,
we call \new{predual pair}
a complete right $\sK$-semimodule $X$ together
with a complete left semimodule $Y$ equipped with a bracket
$\<\cdot,\cdot>$ from $Y\times X$ to a complete $\sK$-bisemimodule
$Z$, such that, for all $x\in X$, the maps
$\sR_x: Y\to Z$, $y\to \<y,x>$
and $\sL_y: X\to Z$, $x\to \<y,x>$
are respectively left and right linear, and continuous.
We shall denote by $(Y,X)$ or simply $Y,X$ this predual pair.
The most familiar choice of $Z$, which corresponds
to ``classical'' bilinear forms, is $Z=\sK$.
The semiring $\sK$ yields another $\sK$-bisemimodule
$Z=\sK\op$, with addition $(x,y)\mapsto \inff\{x,y\}$,
right action $(x,\lambda)\to \lambda \lres  x $,
and left action $(\lambda,x) \to x/\lambda$.

We say that $Y$ \new{separates} $X$
if 
\[ (\forall y\in Y, \<y,x_1>=\<y,x_2>)
\implies x_1=x_2 \enspace,
\]
and that $X$ \new{separates} $Y$
if
\[ (\forall x\in X, \<y_1,x>=\<y_2,x>)
\implies y_1=y_2
\enspace .\]
A predual pair $(Y,X)$
such that $X$ separates $Y$ and $Y$ separates
$X$ is a \new{dual pair}.
This notion is inspired by the dual
pairs which arise in the theory of topological vectors
spaces, 
see~\citename{bourbaki}, Chapter~4,~\citeyear{bourbaki}
or
~\citename{aliprantis}, Chapter~5,~\citeyear{aliprantis}.
\begin{exmp}\label{ex-free2}
The right semimodule $\sK^I$ forms a dual pair with the
left semimodule $\sK^I$ (both were introduced at Example~\ref{ex-free}),
for the canonical bracket $\<a,b>=\supp_{i\in I}a(i)b(i)$.
\end{exmp}
\begin{thm}[Opposite Dual Pair]
Let $X$ denote a complete right $\sK$-semimodule.
Then, the semimodules $X\op,X$ form a dual
pair for the bracket $X\op\times X\to \sK\op$,
$(y,x)\mapsto \<y,x>= x\lres y$.
\end{thm}
\begin{proof}
The bilinearity and continuity of $\<\cdot,\cdot>$
follows from~\eqref{e-inf}, \eqref{e-comp}, \eqref{e-11},
and~\eqref{e-cm}.
Eqn~\eqref{e-separpoints} shows that $X\op$ separates $X$,
and Eqn~\eqref{e-dualseparpoints}
shows that $X$ separates $X\op$.
\end{proof}
A different example of predual pair arises
when considering the (topological) \new{dual} $X'$ of a 
complete semimodule $X$, which is the set of linear
continuous maps $y: X\to \sK$.
The spaces $X',X$ form a predual pair for
the bracket $\<y,x>=y(x)$, and $X$ trivially separates $X'$,
but $X'$ need not separate $X$ (see Example~\ref{ex-nosepar} below).
\begin{exmp}\label{ex-free3}

Consider again the dual pair \((\sK^I,\sK^I)\) of Example~\ref{ex-free2}. 
With any element $a\in \sK^I$ is associated
an element of the dual, $\sL_a: (\sK^I)'$,
$b\mapsto \<a,b>$, and
any element of the dual is of this
form. Thus, $(\sK^I)'$ can be identified to $\sK^I$,
and $(\sK^I)'$ trivially separates $\sK^I$
(indeed, if $b,c\in \sK^I$ are such that
$b(i)\neq c(i)$ for some $i\in I$, the Dirac function 
at point $i$, $\delta_i\in (\sK^I)', \delta_i(d)=d(i)$,
separates $b$ from $c$).
\end{exmp}
\subsection{Involutions}
Given a bracket $\<\cdot,\cdot>$ from
$Y\times X$ to a complete $\sK$-bisemimodule $Z$,
and an arbitrary element $\varphi\in Z$,
we define the maps:
\begin{subequations}
\label{e-def-lrdual}
\begin{align}
X\to Y,\; x\mapsto \ldual{x}&= 
\maxx\set{y\in Y}{\<y,x>\leq \varphi},\\
Y\to X,\;y\mapsto \rdual{y}&= 
\maxx\set{x\in X}{\<y,x>\leq \varphi}
\enspace .
\end{align}
\end{subequations}
Thus, $\ldual{x}=\sR_x\sh(\varphi)$
and $\rdual{y}=\sL_y\sh(\varphi)$.
\begin{prop}
If $(Y,X)$ is a predual pair, then
\begin{subequations}
\label{e-d}
\begin{alignat}{2}
\rdual{(\ldual x)} &\geq x\;, &\qquad \ldual{(\rdual{(\ldual x)})}&= \ldual{x} ,\enspace\forall x\in X\;,
\label{e-d1}\\
\ldual{(\rdual y)}  &\geq y\;, &\qquad 
\rdual{(\ldual{(\rdual{y})})}&= \rdual{y}\; ,\enspace\forall y\in Y \; .
\label{e-d2} 
\end{alignat}
\end{subequations}
\end{prop}
\begin{proof}
We have
\begin{align}
\label{e-dual}
x\leq \rdual{y} \iff \<y,x>\leq \varphi
\iff y\leq \ldual{x}\;  .
\end{align}
Consider now the maps $\sil: Y\to X,\;y\mapsto \rdual{y}$
and $\sir: X\to Y,\;x\mapsto \ldual{x}$. Eqn~\eqref{e-dual}
shows that $\sil: (Y,\leq)\mapsto (X,\leqop)$
is residuated, with $\sil\sh=\sir$. Thus,~\eqref{e-d1}
and~\eqref{e-d2} follow from~\eqref{e-res} and~\eqref{e-invol}.
\end{proof}
We call \new{closed} the elements of $X$ and $Y$
of the form $\rdual{y}$ and $\ldual{x}$, respectively.
We set $\ov X=\set{\rdual y}{y\in Y}$ and $\ov Y=\set{\ldual x}{x\in X}$.
\begin{prop}\label{prop-inf}
The sets of closed elements $\ov X$ and $\ov Y$ 
are complete inf-sub\-semi\-lattices of $X$ and $Y$, respectively,
\end{prop}
\begin{proof}
The set $\ov X$ is the image of the map $\sil:
(Y,\leq )\to (X,\leqop)$ which is residuated,
and, by Lemma~\ref{lem-1}, this image must be
a complete sup-subsemilattice of $X$ for the order $\leqop$,
i.e., a complete inf-subsemilattice of $X$ for the order $\leq$.
\end{proof}
(We warn the reader that the sup laws of $\ov X$
and $\ov Y$ do not coincide with those of $X$ and
$Y$,  in general.)
It follows from~\eqref{e-d} that $\rdual{(\ldual{x})}=x$
(resp.\ $\ldual{(\rdual y)}=y$)
if and only if $x$ (resp.\ $y$) is closed,
hence:
\begin{prop}\label{prop-antitone}
The map $x\mapsto \ldual{x}$ is an anti-isomorphism
of complete lattices $\ov X\to\ov Y$,
with inverse $y\to \rdual{y}$.
\end{prop}
Recall that if $S,T$ are complete lattices, a map $f: S\to T$ is an \new{anti-isomorphism} if, for all $U\subset S$, $f(\supp U)=\inff f(U)$ and $f(\inff U)=\supp f(U)$. A map $f:S\to T$ is 
\new{antitone} if $s\leq s'\implies f(s)\geq f(s')$.
\begin{proof}
We already know that $x\mapsto \ldual{x}$ is an antitone
bijection from $\ov X$ to $\ov Y$ with inverse $y\mapsto \rdual{y}$. A bijective antitone map between complete
lattices whose inverse is antitone
is automatically an anti-isomorphism
of complete lattices. 
\end{proof}

Since $Z$ is a complete $\sK$-bisemimodule,
$\lambda\lres \mu$, and, dually,
$\mu/\nu$ are well defined for $\mu\in Z$ and
$\lambda,\nu \in \sK$.
Considering the predual
pair $(\sK,Z)$
for the bracket $\<\lambda,\mu>=\lambda\mu$
allows us to define $\rdual{\lambda}= \lambda\lres \varphi$.
We define dually $\ldual{\nu}=\varphi/\nu$. 
\begin{prop}
If $x\in X$ and $y\in Y$ are closed, then
\begin{subequations}
\begin{align}
z \lres x &= \rdual{\<\ldual{x},z>} \;,\quad
\forall z\in X\;,\label{e-inva}\\
y / t &= \ldual{\<t,\rdual{y}>}\; ,
\quad\forall t\in Y \;.
\label{e-invb}
\end{align}
\label{e-inv}
\end{subequations}
\end{prop}
\begin{proof}
For all $x\in X$ and
$\nu\in Z$, 
consider the maps $L^{X}_x: \sK\to X,\mu \to x\mu$
and $L^{Z}_\nu : \sK\to Z$, $\mu \to \nu\mu$. We have
$\sL_y\comp L^X_x(\lambda)= \<y,x\lambda>=\<y,x>\lambda$,
for all $\lambda\in \sK$, that is:
\begin{displaymath}
\sL_y \comp L^X_x = L^{Z}_{\<y,x>}\;, \quad \forall y \in Y\;, x\in X \; .
\end{displaymath}
Now, if $x$ is closed, we have $x=\rdual{y}$ for some $y\in Y$,
i.e., $x=\sL_y\sh(\varphi)$. Hence, $z\lres x= z\lres \rdual{y}=
(L^{X}_z)\sh\comp \sL_y\sh(\varphi)=(\sL_y\comp L^X_z)\sh(\varphi)=
(L^{Z}_{\<y,z>})\sh(\varphi)= \rdual{\<y,z>}$, 
which shows~\eqref{e-inva}. We have proved
in passing the following identity, that we tabulate
for further use:
\begin{equation}
\forall y\in Y\;, \enspace z\in X\;,\quad z\lres \rdual{y}= \rdual{\<y, z>}
\; .
\label{e-invc}
\end{equation}
The proof of~\eqref{e-invb} is dual.
\end{proof}
\subsection{Reflexive Semirings}
We say that a complete idempotent semiring $\sK$ 
equipped with a distinguished element $\varphi$
is left (resp.
right) \new{reflexive}
if $\ldual{(\rdual{\lambda})}=\lambda$
(resp.\ $\rdual{(\ldual{\lambda})}=\lambda$),
for all $\lambda\in \sK$, where the operations
$\lambda\mapsto \rdual{\lambda}$,
$\mu\mapsto \ldual{\mu}$ are defined as in~\eqref{e-def-lrdual},
by considering $\sK$ as a bisemimodule over itself,
and taking the bracket $\<\lambda,\mu>=\lambda\mu$.
(The element $\varphi$ need not be unique;
indeed, if $\sK$ is left, or right,
reflexive for $\varphi$,
and if $\lambda$ is invertible, it is
not difficult to check that $\sK$ is
also left (or right) reflexive for $\varphi\lambda$
and $\lambda\varphi$. We shall sometimes write,
more properly, that $(\sK,\varphi)$
is reflexive.)

Using~\eqref{e-injsurj}, together with
$\rdual{\mu}=\sil(\mu)$ and $\ldual{\lambda}= \sir (\lambda)
=\sil\sh(\lambda)$,
we get\begin{subequations}
\begin{align}
\lambda\mapsto \rdual{\lambda}\text{\ is injective}
&\Leftrightarrow \sK \text{\ is left reflexive,}\nonumber\\
&\Leftrightarrow \lambda\mapsto \ldual{\lambda} \text{\ is surjective,}
\label{e-lra}
\\
\lambda\mapsto \ldual{\lambda} \text{\ is injective,}
&\Leftrightarrow \sK \text{\ is right reflexive,}\nonumber\\
&\Leftrightarrow \lambda\mapsto \rdual{\lambda} \text{\ is surjective.}
\label{e-lrb}
\end{align}
\end{subequations}
The interest in reflexive semirings stems
in particular from the following result.
\begin{prop}\label{prop-ref}
If $\sK$ is right reflexive, then the set
of closed elements $\ov X$ is a complete
subsemimodule of $X\op$.
\end{prop}
\begin{proof}
We know from Proposition~\ref{prop-inf} that $\ov X$ is stable by
arbitrary sups for $\leqop$. It remains to check
that for all $x\in\ov X$ and $\lambda\in \sK$,
$\lambda\dotop x=x/\lambda \in \ov X$. By definition
of $\ov X$, we have $x=\rdual{y}=\sL_y\sh(\varphi)$
for some $y\in Y$.
Using~\eqref{e-compose} and the right linearity of
$\<\cdot,\cdot>$, we get $\sL_y\comp R^X_\lambda=R^\sK_\lambda
\comp \sL_y\implies
(R^X_\lambda)\sh \comp \sL_y\sh =
\sL_y\sh \comp (R^\sK_\lambda)\sh$,
hence
$x/\lambda=(R_{\lambda}^X)\sh\comp \sL_y\sh(\varphi)
=\sL_y\sh\comp (R_{\lambda}^\sK)\sh(\varphi)
= \sL_y\sh(\varphi/\lambda)
= \sL_y\sh(\mu\lres \varphi)
$ for some $\mu \in \sK$, since, by~\eqref{e-lrb}, $\mu \mapsto \rdual{\mu}
=\mu\lres \varphi$ is surjective.
Using~\eqref{e-compose}
again,
$x/\lambda= \sL_y\sh \comp (L^{\sK}_\mu)\sh(\varphi)
=(L^{\sK}_\mu\comp \sL_y )\sh(\varphi)
= \sL_{\mu y} \sh(\varphi)=\rdual{(\mu y)}$,
which shows that $x/\lambda\in \ov X$.
\end{proof}
\begin{exmp}
Let us consider once again the dual pair \((\sK^{I},\sK^{I})\)  of Examples~\ref{ex-free}--\ref{ex-free2}--\ref{ex-free3}. Since $\ldual{d}(i)= \varphi/d(i)$,
and since $\rdual{a}(i)= a(i)\lres \varphi$, we see
that all the elements
of the right (resp.~left) semimodule $\sK^I$ are closed as soon as $\sK$ is right (resp.\ left) reflexive.
\end{exmp}
We next exhibit a fundamental class of reflexive idempotent
semirings.
We say that a (non necessarily commutative) semiring
is a \new{semifield} if its non-zero elements
have a multiplicative inverse. A complete idempotent
semiring $\sK$ is never a semifield
(unless $\sK=\{\zero,\unit\}$),
because the maximal element of $\sK$, $\maxxx$,
satisfies $(\maxxx)^2=\maxxx $.
For this reason, we shall call
(in a slightly abusive way) \new{complete semifield}
a complete semiring $\sK$ such that all elements
except $\zero$ and $\maxxx $ have a multiplicative
inverse. For instance,
$\rmaxb=(\R\cup\{\pm\infty\},\max,+)$
is a complete semifield.
\begin{prop}\label{prop-cref}
A complete idempotent semifield $\sK$ is reflexive:
if $\sK=\{\zero,\unit\}$, one must take $\varphi=\zero$,
otherwise, one may take any invertible $\varphi$.
\end{prop}
\begin{proof}
This follows readily from
$\ldual{x}=\varphi x^{-1}$, 
$\rdual{x}=x^{-1}\varphi$, for  $x\not\in \{\zero,\maxxx \}$,
$\ldual{\zero}=\rdual{\zero}=\maxxx $,
$\ldual{(\maxxx)}=\rdual{(\maxxx)}=\zero$. 
\end{proof}
If $G$ is a group,
we denote by $\sK\lser G\rser$
the complete group $\sK$-semialgebra over $G$,
i.e. the free complete $\sK$-semimodule $\sK^G$,
whose elements are denoted as formal sums $\bigoplus_{g\in G} s_g g$
where $\{s_g\}_{g\in G}$ is a family of elements of $\sK$,
equipped with the Cauchy product 
$$(st)_u= \bigoplus_{gh=u\atop g,h\in G} s_g t_h\; .$$
If $\varphi$ is an element of $\sK$, we denote by $\varphi_{\sK\lser G\rser}$
the element of $\sK\lser G\rser$ whose coefficients all are equal to
$\maxxx$,
except the coefficient of the unit, which is equal
to $\varphi_{\sK}$. We also denote by $\varphi_{nn}
\in \sK^{n\times n}$ the matrix
whose diagonal entries are equal to $\varphi$
and whose out-diagonal entries are equal to $\maxxx$.

The abundance of reflexive semirings is shown
by the following immediate property.
\begin{prop}[Transfer Property]
Let $G$ denote a group. If $(\sK,\varphi)$ is a left (or right)
reflexive complete idempotent semiring, 
then so are $(\sK^{n\times n},\varphi_{nn})$ and 
$(\sK\lser G\rser,\varphi_{\sK\lser G\rser})$.\qed
\end{prop}
\begin{prop}\label{prop-closed}
If $\sK$ is reflexive and if $(Y,X)$ form a predual
pair for which $Y$ separates $X$, then, all the elements
of $X$ are closed. 
\end{prop}
\begin{proof}
Since $\sK$ is right reflexive, $\ov X=\set{\rdual{y}}{y\in Y}$
is a complete subsemimodule of $X\op$ (Proposition~\ref{prop-ref}),
hence, applying the Dual Separation Theorem
(Eqn~\eqref{dualsa}) to $V=\ov X\subset X\op$
and to an arbitrary $x\in X\op$, we get, 
\( \forall y\in Y,\quad P\op_{\ov X}(x)\lres \rdual{y}= x\lres \rdual{y}\enspace,
\)
and, using~\eqref{e-invc},
\begin{align}
\forall y\in Y, \quad \rdual{\<y,P\op_{\ov X}(x)>}
= \rdual{\<y,x>} \enspace .
\label{e-u1}
\end{align}
Since $\sK$ is left reflexive, by~\eqref{e-lra},
$\lambda\to \rdual{\lambda}$ is injective,
and, using~\eqref{e-u1}, we get
$\forall y\in Y, \quad \<y,P\op_{\ov X}(x)>
= \<y,x>$. Since $Y$ separates $X$, 
$P\op_{\ov X}(x)=x$, which shows that $x\in \ov X$.
Thus, $X=\ov X$.
\end{proof}
Gathering Proposition~\ref{prop-antitone} and
Proposition~\ref{prop-closed} together with the symmetric
result to Proposition~\ref{prop-closed}, we get:
\begin{cor}\label{c-lovely}
If $(Y,X)$ is a dual pair for a reflexive semiring $\sK$, then 
the map $x\mapsto \ldual{x}$, together with its inverse
$y\mapsto\rdual{y}$,
are anti-isomorphisms of lattices
between $X$ and $Y$.
\end{cor}
\begin{thm}[Hahn-Banach Theorem, Geometric Form]\label{t-hb}
Let $(Y,X)$ denote a predual pair
for a left reflexive semiring $\sK$.
If $V\subset X$ is a complete subsemimodule whose
elements are all closed,
and if $x$ is closed, then, 
\begin{subequations}\label{e-separ2}
\begin{align}
\<\ldual{P_V(x)}, v> = \<\ldual x, v>\;,&\quad\forall v\in V\;,\\\intertext{and}
\<\ldual{P_V(x)}, x> = \<\ldual x, x> &\Leftrightarrow
x\in V \; .
\end{align}
\end{subequations}
\end{thm}
\begin{proof}
Using~\eqref{e-inva}, we rewrite the universal
separation property (Eqn~\ref{e-separ}) as:
\begin{subequations}\label{e-separ3}
\begin{align}
\forall v\in V\;,\quad &\rdual{\<\ldual{P_V(x)}, v>}  = 
\rdual{\<\ldual{x}, v>}\; ,\\\intertext{and}
x\in V\Leftrightarrow\;&
\rdual{\<\ldual{P_V(x)}, x>}  = \rdual{\<\ldual{x}, x>} 
\; .
\end{align}
\end{subequations}
Since $\sK$ is left reflexive, 
as noted in~\eqref{e-lra},
$\lambda\to \rdual{\lambda}$ is injective,
hence,~\eqref{e-separ3} implies~\eqref{e-separ2}.
\end{proof}
A weaker statement,
which is easier to remember,
is the following.
\begin{cor}
If $(Y,X)$ is a predual pair for a reflexive complete semiring $\sK$
such that $Y$ separates $X$,
if $V$ is a complete subsemimodule of $X$, and if $x\in X$,
then, the Hahn-Banach type property~\eqref{e-separ2}
holds.\qed
\end{cor}
\subsection{Representation of Linear Forms}
\label{sec-rep}
We now study the dual pair $(X',X)$.
The following result characterizes the linear
form $\ldual{x}$.
\begin{thm}\label{th-1} 
Let $\sK$ be a complete idempotent reflexive semiring,
let $X$ be a complete $\sK$ semimodule, and consider
the dual pair $(X',X)$ equipped with its canonical bracket.
Then, 
\begin{displaymath}
\ldual{x}(y)= \varphi/(y\lres x)\;,\quad\forall x,y\in X \; .
\end{displaymath}
\end{thm}
\begin{proof}
If $f\in X'$ is such that $f(x)\leq \varphi$, we get
from $x \geq y(y\lres x)$ that $\varphi\geq f(x)\geq f(y)(y \lres x)$,
hence $f(y)\leq \varphi/(y\lres x)$, for all
$y\in X$.
Thus, $\ldual{x}(y)\leq  \varphi/(y\lres x)$.
To show that the equality holds,
it suffices to show that the map $g:X\to \sK$, $y \mapsto \varphi/(y\lres x)$
is linear continuous and satisfies $g(x)\leq \varphi$. 
Since $g(x)=\varphi/(x\lres x)\leq \varphi/\unit=\varphi$, the latter
condition is satisfied. If $\sK$ is reflexive,
the map $\sK\to\sK,\;\lambda\mapsto \varphi/\lambda$,
which is an anti-isomorphism of lattices, 
sends arbitrary infs to arbitrary sups,
and conversely:
\begin{subequations}
\label{eq-ant}
\begin{align}
\varphi/(\inff \Lambda) &=\supp(\varphi/\Lambda)\;,\quad\forall \Lambda\subset \sK \; ,
\label{eq-anti}\\
\varphi/(\supp \Lambda) &=\inff(\varphi/\Lambda)\;,\quad\forall \Lambda\subset \sK \; ,
\label{eq-anti2}
\end{align}
\end{subequations}
(the residuation equality~\eqref{eq-anti2} holds even
if the complete idempotent semiring $\sK$ is not reflexive).
Using~\eqref{eq-anti} and~\eqref{e-11},
we get that for all $V\subset X$,
\[
\varphi/((\supp V) \lres x)=\varphi/(\inff (V\lres x))=\supp (\varphi/(V\lres x)) \;,
\]
which shows that $g$ preserves arbitrary sups. It remains to show
that $g(y\lambda)=g(y)\lambda$, for all $y\in X$, $\lambda \in \sK$.
Since
\[
g(y\lambda)=\varphi/((y\lambda)\lres x)=
\varphi/(\lambda\lres (y\lres x)) \; ,
\]
it suffices to show that $\varphi/(\lambda \lres \alpha)=(\varphi/\alpha)\lambda$ holds for all $\alpha\in \sK$. Since $\sK$ is reflexive,
we can write $\alpha=\beta\lres \varphi$, with
$\beta=\varphi/\alpha$, hence, 
$\varphi/(\lambda \lres \alpha)
= \varphi/(\lambda \lres (\beta\lres \varphi))
= \varphi/((\beta\lambda)\lres \varphi)
= \beta\lambda= (\varphi/\alpha )\lambda $.
\end{proof}
\begin{cor}[$X'$ separates $X$]\label{cor-sep}
If $\sK$ is a complete idempotent reflexive semiring
and if $X$ is a complete $\sK$ semimodule, then $X'$ separates
$X$. 
\end{cor}
\begin{proof}
Let $x,y\in X$. If $f(x)=f(y)$ for all $f\in X'$,
we have in particular, $\ldual{x}(x)=\ldual{x}(y)$.
Since $\lambda \mapsto \varphi/\lambda,\sK\to\sK$ is injective,
we get $x\lres x = y \lres x$,
hence $x\geq y(y\lres x)= y(x\lres x)\geq y$,
which shows that $x\geq y$. By symmetry, $y\geq x$.
\end{proof}
\begin{exmp}\label{ex-nosepar}
The following counterexample shows that, when $\sK$ is not reflexive,
$X'$ need not separate $X$.

Consider the semiring $\nmaxb=\{ \N
\cup\{-\infty,+\infty\},\max,+,0,-\infty\}$ which is complete.
$X=\{\Z\cup\{-\infty,+\infty\},\max\}$ is a complete
$\nmaxb$-semimodule for the action $(x,\lambda) \mapsto x+\lambda$
(with the convention $-\infty+\infty=-\infty$).
Let us prove that $X'$, the set of $\nmaxb$-linear maps from $X$ to
$\nmaxb$, consists only of the two following elements~:
\begin{enumerate}
\item $x\in X \mapsto -\infty$ ;
\item $x\in X \mapsto x+\infty$ .
\end{enumerate}
Let $\phi\in X'$ a linear map and let us assume that it takes only
finite values on $\Z$. Then,
for all $p\in \Z$,
$$\phi(p)=\phi(p-n)+n\geq n\;,\quad\forall n \in \N\;,$$
therefore $\phi(p)\geq \supp n=+\infty$ which is a contradiction.

Let us assume that there exists $p\in\Z$ such that $\phi(p)=-\infty$.
By monotony of $\phi$, $\phi(q)=-\infty$ for all $q\leq p$.
Moreover $\phi(p+n)=\phi(p)+n=-\infty$, for $n\in \N\cup\{+\infty\}$,
which implies $\phi(x)=-\infty$ for all $x$. 

Let us assume that there exists $p\in\Z$ such that $\phi(p)=+\infty$.
 By monotony of $\phi$, $\phi(q)=+\infty$ for all $q\geq p$.
Moreover $\phi(p)=\phi(p-n)+n=+\infty$, for $n\in\N$, which shows that
 $\phi(x)=+\infty$ for all $x\in \Z$, and
 $\phi(-\infty)=-\infty$,
since $\phi$ is linear, so that 
$\phi(x)=x+\infty$ for all $x \in \nmaxb\; .$

Thus, any linear form on $X$ is constant on the set
of finite elements of $X$, which shows that $X'$
does not separate $X$.
\end{exmp}
Since $X'$ separates $X$ and $X$ separates $X'$, we get
as an immediate corollary of Theorem~\ref{th-1}, 
and Corollary~\ref{c-lovely},
the following Riesz representation theorem,
which extends~\cite[Theorem~5.2]{litvinov00}.

\begin{cor}[Riesz Representation Theorem]\label{cor-2}
Let $\sK$ denote a complete idempotent reflexive semiring,
and $X$ a complete $\sK$-semimodule.
Then, any continuous linear form $f\in X'$ can be represented
as 
\begin{equation}
f(y)=\ldual{x}(y)=  \varphi/(y\lres x) \;, \quad \forall y\in X\;,
\label{eq-0}
\end{equation}
for some $x\in X$, and the unique $x\in X$ which satisfies~\eqref{eq-0}
is equal to $\rdual{f}$. 
\end{cor}
We get as a last, immediate corollary, the following
extension of~\cite[Theorem~5.3]{litvinov00}. 
\begin{cor}[Hahn-Banach Theorem, Analytic Form]
\label{cor-hb-ana}
If $\sK$ is a complete idempotent reflexive semiring,
and if $V$ is a complete subsemimodule of a complete
$\sK$-semimodule $X$, then any continuous linear form
defined on $V$ has a continuous extension to $X$.
\end{cor}
\begin{exmp}[Complete Semilattices]\label{ex-boolean}
A complete sup-semilattice $(X,\leq)$ can
be thought of as a complete semimodule over
the Boolean semiring $\B=\{\zero,\unit\}$,
with addition $(x,y)\mapsto\supp\{x,y\}$
and action $x\unit=x$ and $x\zero=\minn X$.
The dual $X'$ is the set
of maps $x':X\to \{\zero,\unit\}$ which preserve
arbitrary sups. Let us take $\varphi=\zero$,
together with the bracket $\<x',x>=x'(x)$
(as noted in Proposition~\ref{prop-cref},
the Boolean semiring has the exceptional
feature of being reflexive for $\varphi=\zero$).
For any $a\in X$, we have
$\ldual{a}=\supp\set{x'\in X'}{x'(a)= \zero}$,
and it is not difficult to see
that $\ldual{a}(x)=\zero$ if $x\leq a $,
and $\ldual{a}(x)=\unit$, otherwise. 
By Corollary~\ref{cor-sep}, $X'$ separates $X$
and, by Corollary~\ref{c-lovely}, 
$x\mapsto \ldual{x}$ establishes an anti-isomorphism
between the lattices $X$ and $X'$.
An equivalent property was already
noticed by
~\citename{wagneur95}~\citeyear{wagneur95}.
\end{exmp}

\subsection{Application: Duality between Row and Column Spaces}

Let $\sK$ denote a complete reflexive semiring,
and let $A\in \sK^{n\times p}$. The free complete 
semimodules $X=\sK^{p\times 1}$ and $Y=\sK^{1\times n}$
form a predual pair for the bracket $\<y,x>=yAx$.
We have $\rdual{y}=\maxx\set{x}{yAx\leq \varphi}
=(yA)\lres \varphi$, and dually ,
$\ldual{x}=\varphi/(Ax)$.
Hence,
\begin{subequations}
\begin{align}
\ov X&= \set{(yA)\lres \varphi}{y\in Y}
\;,\label{e-ovx}\\
\ov Y &= \set{\varphi/(Ax)}{x\in X}
\;.
\end{align}
\end{subequations}
Let $\rowspace(A)=\set{yA}{y\in Y}$
denote the \new{row space}
of $A$, i.e., the left $\sK$-subsemimodule of $\sK^{1\times p}$ generated
by the rows of $A$, and, dually,
let $\colspace(A)=\set{Ax}{x\in X}$
denote the \new{column space} of $A$.
Since $\sK$ is reflexive, the maps $z\mapsto 
z\lres \varphi= (z_i\lres \varphi)_{1\leq i\leq p}$
and $z\mapsto \varphi/z= (\varphi/z_i)_{1\leq i\leq p}$
are mutually inverse antitone bijections
between $\sK^{1\times p}$ and $\sK^{p\times 1}=X$.
By~\eqref{e-ovx}, $z\mapsto z\lres \varphi$
sends $\rowspace(A)$ to $\ov{X}$, hence,
$\rowspace(A)$ and $\ov{X}$ are
anti-isomorphic lattices. Dually,
$\colspace(A)$ and $\ov{Y}$ are
anti-isomorphic lattices.
By Proposition~\ref{prop-antitone},  
$\ov{X}$ and $\ov{Y}$ are
anti-isomorphic lattices.
Composing anti-isomorphisms, we see that
the map:
\[
\rowspace(A)\to \colspace(A)\;,\quad
z\mapsto \big[\varphi/ \big(A (z\lres\varphi)\big)\big]\lres \varphi
= A (z\lres\varphi)
\]
is an anti-isomorphism of lattices.
We have proved the following
result, which extends a theorem
of Markowsky (see~\cite[Theorem~1.2.3]{kim82}) for Boolean
matrices.
\begin{thm}
The row space and column space
of a matrix with entries in a complete
idempotent reflexive semiring are anti-isomorphic lattices.\qed
\end{thm}
\def\cprime{$'$}

\end{document}